\documentclass[preprint]{elsarticle}
\usepackage{amsmath,amssymb,enumerate,amscd}
\input amssym.def
\usepackage[utf8]{inputenx}
\newcommand{\CP}[1]{\makebox{${\rm C}_p({#1})$}}

\newcommand{\C}[1]{\makebox{${\rm C}({#1})$}}

\newcommand{\USCB}[1]{\makebox{${\rm USC}^*({#1})$}}
\newcommand{\USC}[1]{\makebox{${\rm USC}({#1})$}}

\newcommand{\seq}[2]{\left\langle #1:#2\in\omega\right\rangle}

\newcommand{\seqc}[1]{\makebox{$[#1]_{\scriptstyle seq}$}}

\newcommand{\nula}{\makebox{${\mathbf 0}$}}

\newcommand{\snula}{\makebox{${\scriptstyle\mathbf 0}$}}
\newcommand{\RR}{{\mathbb R}}
\newcommand{\QQ}{{\mathbb Q}}

\newcommand{\constant}[1]{\makebox{${\mathbf #1}$}}

\newtheorem{theorem}{Theorem}
\newtheorem{corollary}[theorem]{Corollary}
\newtheorem{lemma}[theorem]{Lemma}

\newcommand{\pf}{\noindent{\bf Proof. }}

\newcommand{\emm}{\bf}

\begin{document}
\begin{frontmatter}
\title{Spaces of Real Functions, Covers and Dense Subsets}
\author{Lev Bukovsk\'y}
\ead{lev.bukovsky@upjs.sk}
\address{Institute of Mathematics, Faculty of Science,
P.J. \v Saf\'arik University, Jesenn\'a~5, 040~01~Ko\v sice,
Slovakia}
\author{Alexander V. Osipov}
\ead{OAB@list.ru}
\address{Krasovskii Institute of Mathematics and Mechanics, Ural Federal University, \\ Ural State University of Economics, 620219, Yekaterinburg, Russia}
\begin{abstract}
In paper we study relationships between covering properties of a
topological space $X$ and the space
$(USC^*(X),\tau_{\mathcal{B}})$ of bounded upper semicontinuous
functions on $X$ with the topology $\tau_{\mathcal{B}}$ defined by
the bornology $\mathcal{B}$ on $X$. We present characterizations
of various local properties of $(USC^*(X),\tau_{\mathcal{B}})$ in
terms of selection principles related to bornological covers of
$X$. We also investigate topological properties of sequences of
upper dense and upper sequentially dense subsets of
$(USC^*(X),\tau_{\mathcal{B}})$.

\end{abstract}
\end{frontmatter}

\section{Introduction}

In this paper we tries to give results about the relationships of
the properties of some covers of a topological space $X$ and the
properties of the families of real functions on $X$. We tried to
extend the known results also for covers respecting a bornology on
$X$ and spaces of bounded upper semicontinuous functions with the
corresponding topology defined by the bornology. This paper is a
natural generalization the research done in papers
\cite{BL1,BL2,BH,BO,Osi,Osi3,Osi4,Osi6,Osi7} for spaces of
continuous real-valued functions with the topology of pointwise
convergence and with the compact-open topology.

\section{Covers and bornology}

A topological space $(X,\tau)$ is always an infinite Hausdorff
topological space, $\tau$ is the family of open subsets of $X$.
Unexplained notions and terminology are those of R. Engelking
\cite{Eng}.

A family $\mathcal{U}$ of subsets of $X$ is a {\it cover} of $X$
if $\bigcup \mathcal{U}=X$. For some technical reasons, a cover
will be called also {\it o-cover}. A cover $\mathcal{V}\subseteq
\mathcal{U}$ is said to be a subcover of $\mathcal{U}$. If we deal
with a countable cover of $X$, we can consider it a sequence of
subsets. A cover is {\it open} if every element of the cover is an
open set.

We say that a family $\mathcal{V}\subseteq \mathcal{P}(X)$ is a
{\it refinement} of the family $\mathcal{U}\subseteq
\mathcal{P}(X)$ if

$(\forall V\in\mathcal{V})(\exists U\in\mathcal{U}) V\subseteq U$.

A {\it bornology} $\mathcal{B}$ on a topological space $X$ is a
proper ideal of subsets of $X$ such that $\bigcup \mathcal{B}=X$.
A subset $\mathcal{B}_0\subseteq \mathcal{B}$ is a {\it base} of
the bornology $\mathcal{B}$ if for every $B\in \mathcal{B}$ there
exists a $B_0\in \mathcal{B}_0$ such that $B\subseteq B_0$. Note
that a bornology $\mathcal{B}$ has a closed base if and only if
for every $B\in\mathcal{B}$ also $\overline{B}\in \mathcal{B}$.
The smallest bornology on $X$ is the ideal $Fin=Fin(X)$ of all
finite subsets of $X$.

We shall use the following convention. If the lower case letters
$\varphi$ or $\psi$ denote one of the symbols $o$, $\lambda$,
$\omega$ or $\gamma$, then the capital letters $\Phi$ or $\Psi$
denote the corresponding symbol $\mathcal{O}$, $\Lambda$, $\Omega$
of $\Gamma$, respectively, and vice versa.

Let $\mathcal{B}$ be a bornology on a topological space $X$. We
shall consider covers respecting this bornology. We assume that a
{\it bornological cover}, briefly $\mathcal{B}$-$o$-cover, is
identical with an $o$-cover. Similarly, a {\it large bornological
cover} $\mathcal{U}$, briefly $\mathcal{B}$-$\lambda$-cover, is
simply a large cover, i.e., for every $x\in X$ the set $\{U\in
\mathcal{U}: x\in U\}$ is infinite. A cover $\mathcal{U}$ is a
{\it bornological $\omega$-cover}, briefly
$\mathcal{B}$-$\omega$-{\it cover}, if $X\not\in \mathcal{U}$ and
for every $B\in \mathcal{B}$ there exists $U\in \mathcal{U}$ such
that $B\subseteq U$. A cover $\mathcal{U}$ is a {\it bornological
$\gamma$-cover}, briefly $\mathcal{B}$-$\gamma$-cover, if
$\mathcal{U}$ is infinite and for every $B\in \mathcal{B}$ the set
$\{U\in \mathcal{U}: B\not\subseteq U\}$ is finite. If
$\mathcal{U}$ is a $\mathcal{B}$-$\gamma$-cover, then
$\mathcal{U}\setminus \{X\}$ is a $\mathcal{B}$-$\gamma$-cover as
well. So, we can assume that $X$ does not belong to a
$\mathcal{B}$-$\gamma$-cover. We denote by $\Phi_{\mathcal{B}}(X)$
the family of all open $\mathcal{B}$-$\varphi$-covers of $X$ for
$\varphi=o,\lambda,\omega,\gamma$.

If $\mathcal{B}=Fin(X)$, then a $\mathcal{B}$-$\varphi$-cover is
the classical $\varphi$-cover and

$\Gamma_{Fin}(X)=\Gamma(X)$, $\Omega_{Fin}(X)=\Omega(X)$.

Evidently,  $\Gamma_{\mathcal{B}}(X)\subseteq
\Omega_{\mathcal{B}}(X)\subseteq
\Lambda_{\mathcal{B}}(X)=\Lambda(X)\subseteq
\mathcal{O}_{\mathcal{B}}(X)=\mathcal{O}(X)$.

Let the family $\mathcal{V}\subseteq \mathcal{P}(X)$ be a
refinement of the family $\mathcal{U}\subseteq \mathcal{P}(X)$. If
$\mathcal{V}$ is an $\mathcal{B}$-$o$- or an
$\mathcal{B}$-$\omega$-cover, then $\mathcal{U}$ is such a cover
as well. This is not true for $\mathcal{B}$-$\lambda$- and an
$\mathcal{B}$-$\gamma$-covers. If we add finitely many subsets of
$X$ to a $\mathcal{B}$-$\gamma$-cover, we obtain a
$\mathcal{B}$-$\gamma$-cover. Moreover, each infinite subset of a
$\mathcal{B}$-$\gamma$-cover is a $\mathcal{B}$-$\gamma$-cover as
well. Omitting finitely many elements of an
$\mathcal{B}$-$\lambda$- or an $\mathcal{B}$-$\omega$-cover, we
obtain a cover of same type. This is not true for
$\mathcal{B}$-$o$-cover.

A $\mathcal{B}$-$\varphi$-cover $\mathcal{U}$ is {\it shrinkable}
if there exists an open $\mathcal{B}$-$\varphi$-cover
$\mathcal{V}$ such that

$(\forall V\in \mathcal{V})(\exists U_{V}\in \mathcal{U}\setminus
\{X\})\overline{V}\subseteq U_{V}$.

The family $\{U_V: V\in \mathcal{V}\}\subseteq \mathcal{U}$ is a
$\mathcal{B}$-$\varphi$-cover as well. The family of all open
shrinkable $\mathcal{B}$-$\varphi$-covers of $X$ will be denoted
by $\Phi^{sh}_{\mathcal{B}}(X)$, or simply
$\Phi^{sh}_{\mathcal{B}}$.

Similarly to F. Gerlits and Z. Nagy \cite{GN}, we define: $X$ has
{\it the property $(\epsilon_{\mathcal{B}})$} if every open
$\mathcal{B}$-$\omega$-cover contains a countable
$\mathcal{B}$-$\omega$-subcover.

G. Beer and S. Levy in \cite{BL1} introduced the notion of a
strong $\mathcal{B}$-cover of a metric space. It is easy to define
that notion for a uniform space. So, let $(X,\mathcal{W})$ be a
uniform space. The ball about $B\subseteq X$ and radius $V\in
\mathcal{W}$ is the set

$\mathbb{B}(B,V)=\{x\in X: (\exists y\in B)(x,y)\in V\}$.

If $B=\{x\}$, we write simply $\mathbb{B}(x,V)$.

Let $\mathcal{B}$ be a bornology on $X$. An open cover
$\mathcal{U}$ is a {\it strong $\mathcal{B}$-$\omega$-cover},
briefly a {\it $\mathcal{B}$-$\omega^s$-cover}, if $X\not\in
\mathcal{U}$ and for every $B\in \mathcal{B}$ there exists a
$U\in\mathcal{U}$ and a $V\in \mathcal{W}$ such that
$\mathbb{B}(B,V)\subseteq U$. An open cover $\mathcal{U}$ is a
{\it strong $\mathcal{B}$-$\gamma$-cover}, briefly a {\it
$\mathcal{B}$-$\gamma^s$-cover}, if $\mathcal{U}$ is infinite and
for every $B\in \mathcal{B}$ the set $\{U\in\mathcal{U}: \neg
(\exists V\in \mathcal{W})\mathbb{B}(B,V)\subseteq U\}$ is finite.
As above, we can assume that $X$ does not belong to a
$\mathcal{B}$-$\gamma^s$-cover. We denote by
$\Omega_{\mathcal{B}}^s(X)$ and $\Gamma_{\mathcal{B}}^s(X)$ the
family of all open $\mathcal{B}$-$\omega^s$-covers and open
$\mathcal{B}$-$\gamma^s$-covers of $X$, respectively. Then, we
have $\Gamma_{\mathcal{B}}^s(X)\subseteq \Omega_{\mathcal{B}}^s(X)
\subseteq \mathcal{O}(X)$.

Similar as above, we have $\Gamma_{Fin}^s(X)=\Gamma(X)$,
$\Omega_{Fin}^s(X)=\Omega(X)$.

One can easily see that for $\Phi=\Omega, \Gamma$ we have
$\Phi_{\mathcal{B}}^s(X)\subseteq \Phi_{\mathcal{B}}(X)\subseteq
\Phi(X)$.

Both types of covers suggest to introduce corresponding topology
on $^X\mathbb{R}$.

The topology $\tau_{\mathcal{B}}$ is defined by typical
neighborhoods of a function $h\in ^X\mathbb{R}$ of the form

\begin{equation}\label{neig}
\mathcal{N}_{\mathcal{B},\epsilon}(h)=\{f\in ^X\mathbb{R}:
(\forall x\in B)|h(x)-f(x)|<\epsilon\}
\end{equation}

for a set $B\in \mathcal{B}$ and $\epsilon>0$.

The product topology $\tau_p$ on $^X\mathbb{R}$ is actually the
topology $\tau_{Fin}$.

For a uniform space $(X,\mathcal{W})$, the topology related to
$\mathcal{B}$-$\varphi^s$-covers $\tau_{\mathcal{B}}^s$ defined by
typical neighborhoods of a function $h\in ^X\mathbb{R}$ of the
form

\begin{equation}\label{neig1}
\mathcal{N}^s_{\mathcal{B},\epsilon}(h)=\{f\in ^X\mathbb{R}:
(\exists V\in \mathcal{W})(\forall x\in
\mathbb{B}(B,V))|h(x)-f(x)|<\epsilon\}
\end{equation}
 for $B\in \mathcal{B}$ and $\epsilon>0$.

One can easily see that $\tau_p\subseteq
\tau_{\mathcal{B}}\subseteq \tau_{\mathcal{B}}^s$.

\section{Families of real functions}

Let $\mathcal{B}$ be a bornology on $X$. Similarly as
in~\cite{BL2} we introduce the following properties of a~family
$F$ of real functions and a~function $h\in {}^X\RR$:
\[
\begin{array}{ll}
({\mathcal O}_h)_{\mathcal{B}}&h(x)\in\overline{\{f(x):f\in F\}}\textnormal{\ for\ every\ } x\in X.\\
(\Omega_h)_{\mathcal{B}}&h\notin F\textnormal{\  and\ }h\in\overline{F} \ in\ the\ topology\ of \tau_{\mathcal{B}}.\\
(\Gamma_h)_{\mathcal{B}}&F\textnormal{\ is\ infinite and for\ every\ }\varepsilon>0\textnormal{\ and for every\ }B\in \mathcal{B}\\
{}&\textnormal{the\ set\ }\{f\in F:(\exists x\in B)\vert
f(x)-h(x)\vert\geq \varepsilon\}\textnormal{\ is\ finite}.
\end{array}
\]

Omitting $h$ from a set $F$ with $(\Gamma_h)_{\mathcal{B}}$, we
obtain $(\Gamma_h)_{\mathcal{B}}\rightarrow
(\Omega_h)_{\mathcal{B}}\rightarrow ({\mathcal
O}_h)_{\mathcal{B}}$.

One can easily see that for $\Phi=\mathcal{O},\Omega,\Gamma$ we
have

If $\langle F, +\rangle$, $F\subseteq ^X\mathbb{R}$ is a group,
then $F$ has the property $(\Phi_h)_\mathcal{B}$ if and only if
$F$ has the property $(\Phi_{h+f})_\mathcal{B}$ for every $f\in
F$.

\medskip
The set ${}^X\RR$ of all real function defined on $X$ is endowed
with the product topology. Thus, a~typical neighborhood of
a~function \makebox{$g\in{}^X\RR$} is the set
\begin{equation}\label{neig2}
V=\{h\in{}^X\RR:\vert h(x_j)-g(x_j)\vert<\varepsilon:j=0,\dots,k\}
\end{equation}
where $\varepsilon$ is a~positive real and $x_0,\dots,x_k\in X$.
A~sequence of real functions $\seq{f_n}{n}$ converges to a~real
function $f$ in this topology if it converges pointwise, i.e., if
$\lim_{n\to\infty}f_n(x)=f(x)$ for each $x\in X$.
\par
Similarly as in~\cite{BL2} we introduce the following properties
of a~family $F$ of real functions and a~function $h\in {}^X\RR$:
\[
\begin{array}{ll}
({\mathcal O}_h)&h(x)\in\overline{\{f(x):f\in F\}}\textnormal{\ for\ every\ } x\in X.\\
(\Omega_h)&h\notin F\textnormal{\  and\ }h\in\overline{F} \ in\ the\ topology\ of\ {}^X\RR.\\
(\Gamma_h)&F\textnormal{\ is\ infinite and for\ every\ }\varepsilon>0\textnormal{\ and for every\ }x\in X\\
{}&\textnormal{the\ set\ }\{f\in F:\vert f(x)-h(x)\vert\geq
\varepsilon\}\textnormal{\ is\ finite}.
\end{array}
\]


Let $H\subseteq {}^X\RR$. We set
\[\Phi_h(H)=\{F\subseteq H:F\textnormal{\ possesses\ }(\Phi_h)\land (\forall f\in F)\,(f\geq h\land f-h\in H)\}.\]
\par
For a~real $a$, we denote by $\constant{a}$ the constant function
on $X$ with value~$a$. For simplicity for $g\in {}^X\RR$, instead
of $g+\constant{a}$ or $g-\constant{a}$ we shall write $g+a$ or
$g-a$, respectively. Similarly for $\min\{\constant{a},g\}$ or
$\max\{\constant{a},g\}$. If $F\subseteq {}^X\RR$, then
\[F^+=\{f\in F:f\geq 0\},\ \ F^*=\{f\in F :f\textnormal{\ is bounded}\}.\]
\par
$C(X)$ or $\USC{X}$ denote the set of all real continuous or upper
semicontinuous functions\footnote{A function $f:X\longrightarrow
\RR$ is said to be upper semicontinuous if for every real~$a$ the
set $\{x\in X:f(x)<a\}$ is open.}  defined on the~topological
space $X$. Instead of $C(X)^*$ or $\USC{X}^*$ we write $C^*(X)$ or
$\USCB{X}$, respectively.

A~set $F\subseteq H\subseteq \mathbb{R}^X$ is {\emm sequentially
dense in} $H$ if $H\subseteq \seqc{F}$. The~set $F\subseteq
H\subseteq ^X\mathbb{R}$ is {\emm countably dense in} $H$ if for
every function~$f\in H$ there exists a~countable set $G\subseteq
F$ such that $f\in\overline{G}$. As~obviously, the~set $F$ is
{\emm dense in} $H$  if $H\subseteq \overline{F}$. Finally,
the~set $F$ is {\emm pointwise dense in} $A\subseteq\RR$ if
$A\subseteq \overline{\{f(x):f\in F\}}$ for each $x\in X$ ({\bf
$1$-dense set} in terminology of \cite{Osi2,Osi3}). We set
\[
\begin{array}{lll}
{}&{\mathcal S}(H)=&\!\!\!\! \{F\subseteq H:F\textnormal{\ is sequentially dense in\ }H\}, \nonumber \\
{}&{\mathcal D}(H)=&\!\!\!\! \{F\subseteq H:F\textnormal{\ is dense in\ }H\} , \nonumber \\
{}&{\mathcal P}(H)=&\!\!\!\! \{F\subseteq H:F\textnormal{\ is
pointwiae dense in\ }H\} . \nonumber
\end{array}
\]
Then
\[{\mathcal S}(H)\subseteq {\mathcal D}(H)\subseteq {\mathcal P}(H).\]
\par
Evidently a~sequentially dense set is countably dense as well. By
Tong Theorem, see, e.g., \cite{Eng}, if $X$ is perfectly normal
topological space then every (bounded) upper semicontinuous
function is a~limit of a~non-increasing sequence of (bounded)
continuous functions. Thus for a~perfectly normal topological
space $X$  the set $(C(X),\tau_p)$ is sequentially dense in
$(USC(X),\tau_p)$. Then the set $(C^*(X),\tau_p)$ is sequentially
dense in $(\USCB{X},\tau_p)$ as well.

We shall modify the results in~\cite{BO} for bornological covers
and topologies. Actually, we follow the paper~\cite{BO}.
\par
Note that the following notions depend on the topology we consider on $H$. For our purpose we assume that $H$ is endowed with $\tau_{\mathcal B}$ topology. A~set $F\subseteq H\subseteq {}^X\RR$ is {\emm upper sequentially dense in\ }$H$ if for every $f\in H$ there exists a~sequence $\seq{h_n}{n}$ of elements of $F$ such that $h_n\to f$, $h_n\geq f$ and $h_n-f\in H$ for each $n\in\omega$.  A~set $F\subseteq H$ is {\emm upper dense in\ }$H$ if for every $f\in H$ the set $\{h\in F:h\geq f\land h-f\in H\}$ is dense in the set~$\{h\in H:h\geq f\}$. Similarly we define the notions "{\emm upper sequentially dense in\ }$\langle H,\tau_{\mathcal B}\rangle$", "{\emm upper sequentially dense in\ }$\langle H,\tau_{{\mathcal B}^s}\rangle$", "{\emm upper dense in\ }$\langle H,\tau_{\mathcal B}\rangle$", "{\emm upper dense in\ }$\langle H,\tau_{{\mathcal B}^s}\rangle$".
\par
One can easily see that if a~set $F\subseteq \USCB{X}$ is upper
dense in $\USCB{X}$, then for every continuous function $f$ the
set of upper semicontinuous functions \mbox{$\{h-f: h\in F\land
h\geq f\}$} is upper dense in $\USCB{X}^+$. If the~set $F$ is
upper sequentially dense in $\USCB{X}$, then for every \mbox{$f\in
C(X)$} the set $\{h-f:h\in F\land h\geq f\}$ is upper sequentially
dense in $\USCB{X}^+$. We set
\[
\begin{array}{lll}
{}&{\mathcal S}^{\uparrow}_{\mathcal B}(H)=&\{F\subseteq H:F\textnormal{\ is upper sequentially dense in\ }\langle H,\tau_{\mathcal B}\rangle\}, \nonumber \\
{}&{\mathcal D}^{\uparrow}_{\mathcal B}(H)=&\{F\subseteq H:F\textnormal{\ is upper dense in\ }\langle H,\tau_{\mathcal B}\rangle\} , \nonumber \\
{}&{\mathcal P}^{\uparrow}_{\mathcal B}(H)=&{\mathcal P}(H)\nonumber
\end{array}
\]
and similarly for $\tau^s$.
Then
\[{\mathcal S}^{\uparrow}_{\mathcal B}(H)\subseteq {\mathcal D}^{\uparrow}_{\mathcal B}(H)\subseteq {\mathcal P}^{\uparrow}_{\mathcal B}(H).\]
\par
If ${\mathcal B}={\rm Fin}$ then we simply omit the index and we obtain the notation of~\cite{BO}.
\par
We introduce the following notations. If $\Phi=\Gamma$ then $\tilde\Phi={\mathcal S}$. If $\Phi=\Omega$ then $\tilde\Phi={\mathcal D}$ and if $\Phi={\mathcal O}$ then $\tilde\Phi={\mathcal P}$. Similarly for~$\tilde\Phi^{\uparrow}$.
\par
Note that by definitions we have immediately
\begin{equation}\label{PhiPhih}
(\forall F\in \tilde\Phi_{\mathcal B}(H)^{\uparrow})(\forall h\in H)(\exists G\subseteq F)\,G\in\Phi_{h,{\mathcal B}}(H).
\end{equation}

\section{Dense selectors of $\langle USC^*(X),\tau_{\mathcal{B}}\rangle$ and $\langle
USC^*(X),\tau_{\mathcal{B}^s}\rangle$}
\par

In \cite{BO} the authors introduced the following set of real functions:
\begin{equation}\label{SU}
S({\mathcal U})=\{f_{U,g}:U\in{\mathcal U}\land g\in\USCB{X}\},
\end{equation}
where
\begin{equation}\label{fUh}
f_{U,g}(x)=\left\{
\begin{array}{ll}
g(x)&\mbox{\ if\ } x\in U,\\
g(x)+1+\sup \vert g\vert&\mbox{\ otherwise}.
\end{array}
\right.
\end{equation}
We show the basic properties of the families $S({\mathcal U})$.
\begin{lemma}\label{lem1}\ \ \
\begin{enumerate}
\item[{\rm a)}] If ${\mathcal U}$ is an~open ${\mathcal B}$-$\omega$-cover, then the family $S({\mathcal U})\subseteq\USC{X}$ is upper dense in $\langle \USCB{X},\tau_{\mathcal B}\rangle$.
\item[{\rm b)}] If ${\mathcal U}$ is an~open ${\mathcal B}$-$\gamma$-cover, then the family $S({\mathcal U})\subseteq\USC{X}$ is upper sequentially dense in $\langle \USCB{X},\tau_{\mathcal B}\rangle$.
\end{enumerate}
If $X$ is a~metric space then
\begin{enumerate}
\item[{\rm c)}] If ${\mathcal U}$ is an~open ${\mathcal B}^s$-$\omega$-cover, then the family $S({\mathcal U})\subseteq\USC{X}$ is upper dense in $\langle \USCB{X},{\tau_{\mathcal B}^s}\rangle$.
\item[{\rm d)}] If ${\mathcal U}$ is an~open ${\mathcal B}^s$-$\gamma$-cover, then the family $S({\mathcal U})\subseteq\USC{X}$ is upper sequentially dense in $\langle \USCB{X},\tau_{{\mathcal B}^s}\rangle$.
\end{enumerate}
\end{lemma}
\pf
We prove only parts a) and d), since the others can be proved in a~very similar way.
\par
One can easily see that $f_{U,g}\geq g$  and $f_{U,g}$ is bounded upper
semicontinuous for $g\in \USCB{X}$.
\par
We show that if ${\mathcal U}$ is an~open ${\mathcal
B}$-$\omega$-cover then $S({\mathcal U})$ is upper dense in
$\langle \USCB{X},\tau_{\mathcal B}\rangle$. Assume that
$g\in\USCB{X}$. If $[B,\varepsilon](g)$ is a~neighborhood of $g$
defined by \eqref{neig}, then there exists a~$U\in {\mathcal U}$
such that $B\subseteq U$. Then $f_{U,g}(x)=g(x)$ for $x\in B$.
Hence $f_{U,g}\in [B,\varepsilon](g)$, $f_{U,g}\geq g$ and
$f_{U,g}-g\in\USCB{X}$. Thus $S({\mathcal U})$ is upper dense in
$\langle \USCB{X},\tau_{\mathcal B}\rangle$.
\par
Now we show d), i.e.,  if ${\mathcal U}$ is an~open ${\mathcal
B}^s$-$\gamma$-cover then $S({\mathcal U})$ is upper sequentially
dense in $\langle \USCB{X},\tau_{{\mathcal B}^s}\rangle$. Indeed,
let $g\in\USCB{X}$. Let $\{U_i:i\in\omega\}$ be a~countable
$\gamma$-subcover of ${\mathcal U}$. For $i\in\omega$, we let
$g_i=f_{U_i,g}\in S({\mathcal U})$. We show that the sequence
$\seq{g_i}{i}$ converges to $g$. Let $V$ be a~neighborhood of $g$
defined by \eqref{neig1}. Since $\{U_i:i\in\omega\}$ is
a~${\mathcal B}^s$-$\gamma$-cover, there exists an~$i_0$ such that
$B\subseteq U_i$ for $i\geq i_0$. Then we have $g_i(x)=g(x)$ for
$x\in B$ and $i\geq i_0$. Therefore the elements of the sequence
$\seq{g_i}{i}$ belong to $V$ for~$i\geq i_0$. As above, $g_i\geq
g$ and $g_i-g\in\USCB{X}$. Thus $S({\mathcal U})$ is upper
sequentially dense in $\langle \USCB{X},\tau_{{\mathcal
B}^s}\rangle$. \qed
\begin{theorem}\label{PhiD}
Let $\Phi=\Omega,\Gamma$. Then the following are equivalent:
\begin{enumerate}
\item[{\rm a)}] $\langle \USCB{X},\tau_{\mathcal B}\rangle$
satisfies the selection principle {\rm S}${}_1(\tilde
\Phi^{\uparrow}_{\mathcal B},{\mathcal D})$. \item[{\rm b)}]
$\langle \USCB{X},\tau_{\mathcal B}\rangle$ is separable and the
topological space~$X$ possesses the covering property {\rm
S}${}_1(\Phi_{\mathcal B},\Omega_{\mathcal B})$. \item[{\rm c)}]
$\langle \USCB{X},\tau_{\mathcal B}\rangle$ is separable and
satisfies the selection principle
\newline
{\rm S}${}_1(\tilde\Phi_{h,{\mathcal B}},\Omega_{h,{\mathcal B}})$
for every $h\in \USCB{X}$. \item[{\rm d)}] $\langle
\USCB{X},\tau_{\mathcal B}\rangle$ is separable and satisfies the
selection principle
\newline
{\rm S}${}_1(\tilde \Phi^{\uparrow}_{\mathcal
B},\Omega_{h,{\mathcal B}})$ for every $h\in \USCB{X}$.
\end{enumerate}
\end{theorem}
\pf
\newline
${\rm a)}\to {\rm b)}$. Let $\{{\mathcal
U}_n:n\in\omega\}\subseteq \Phi_{\mathcal B}$. We may assume that
${\mathcal U}_{n+1}$ is a~refinement of ${\mathcal U}_n$ for each
$n\in\omega$. If $\Phi=\Gamma$ we may also assume that  for every
$n\in\omega$, the cover ${\mathcal U}_n$ is a~countable family
$\{U_i^n:i\in\omega\}$.
\par
 For every $n\in\omega$ we set
\begin{equation}\label{Sn}
S_n=S({\mathcal U}_n).
\end{equation}
By Lemma \ref{lem1} we have $S_n\in\tilde\Phi^{\uparrow}_{\mathcal
B}$. Thus, by the selection principle {\rm S}${}_1(\tilde
\Phi^{\uparrow}_{\mathcal B},{\mathcal D})$, for every
$n\in\omega$ we obtain an~$f_{U_n,h_n}\in S_n$ such that
$\{f_{U_n,h_n}:n\in\omega\}$ is dense in $\langle
\USCB{X},\tau_{\mathcal B}\rangle$. We show that
$\{U_n:n\in\omega\}$ is a ${\mathcal B}$-$\omega$-cover.
\par
Let $B\in {\mathcal B}$. Consider the open non-empty set
\[U=\{g\in USC^*(X):\vert g(x)\vert<1/2\textnormal{\ for\ }x\in B\}\]
Since the set  $\{f_{U_n,h_n}:n\in\omega\}$ is dense in $\langle
USC^*(X), \tau_{\mathcal B}\rangle$, there exists an~$n$ such that
$f_{U_n,h_n}\in U$. Since $\vert f_{U_n,h_n}(x)\vert <1/2$ for
$x\in B$, by \eqref{fUh} we obtain $B\subseteq U_n$.
\par
The implication ${\rm b)}\to{\rm c)}$ follows by Corollary.
\par
The implication ${\rm c)}\to{\rm d)}$ is obvious by \eqref{PhiPhih}.
\newline
${\rm d)}\to{\rm a)}$. We assume that $\langle
\USCB{X},\tau_{\mathcal B}\rangle$ is separable and satisfies the
selection principle {\rm S}${}_1(\tilde \Phi^{\uparrow}_{\mathcal
B},\Omega_{h,{\mathcal B}})$ for every $h\in USC^*(X)$. Thus,
there exists a~ countable set $D=\{d_n:n\in\omega\}$ dense in
$\langle \USCB{X},\tau_{\mathcal B}\rangle$.
 Let $\{S_{n,m}:n,m\in\omega\}$ be a~sequence of
subsets of $\langle \USCB{X},\tau_{\mathcal B}\rangle$ such that
$S_{n,m}\in \tilde \Phi^{\uparrow}_{\mathcal B}$ for each $n,m\in
\omega$. For every $n\in\omega$ we apply the sequence selection
principle {\rm S}${}_1(\tilde \Phi^{\uparrow}_{\mathcal
B},\Omega_{{d_n},{\mathcal B}})$ to the sequence
$\seq{S_{n,m}}{m}$ and for every $m\in\omega$ we obtain
$d_{n,m}\in S_{n,m}$ such that
$d_n\in\overline{\{d_{n,m}:m\in\omega\}}$. Then
$\{d_{n,m}:n,m\in\omega\}$ is  dense in $\langle
\USCB{X},\tau_{\mathcal B}\rangle$. \qed
\par

\medskip

Analogously to the proof of Theorem \ref{PhiD} we get the
following theorem.
\begin{theorem}\label{PhiD1}
Let $\Phi=\Omega$ or $\Phi=\Gamma$. Assume that $C_{\mathcal{
B}}(X)$ is countably dense in $\langle \USCB{X},\tau_{\mathcal
B}\rangle$. Then for any couple $\langle\Phi,\Psi\rangle$
different from $\langle \Omega_{\mathcal B},{\mathcal O}\rangle$,
the following are equivalent:
\begin{enumerate}
\item[{\rm a)}] $\langle \USCB{X},\tau_{\mathcal B}\rangle$
satisfies the selection principle {\rm S}${}_1(\tilde
\Phi^{\uparrow}_{\mathcal B},{\mathcal D})$, \item[{\rm b)}]
$\langle \USCB{X},\tau_{\mathcal B}\rangle$ is separable  and the
topological space~$X$ possesses the covering property {\rm
S}${}_1(\Phi_{\mathcal B},\Omega_{\mathcal B})$. \item[{\rm c)}]
$\langle \USCB{X},\tau_{\mathcal B}\rangle$ is separable and
satisfies the selection principle
\newline
{\rm S}${}_1(\Phi_{{\snula},{\mathcal
B}},\Omega_{{\snula},{\mathcal B}})$. \item[{\rm d)}] $\langle
\USCB{X},\tau_{\mathcal B}\rangle$ is separable and satisfies the
selection principle
\newline
{\rm S}${}_1(\tilde \Phi^{\uparrow}_{\mathcal
B},\Omega_{{\snula},{\mathcal B}})$.
\end{enumerate}
\end{theorem}
\pf
\par
We prove only the implication ${\rm d)}\to{\rm a)}$. The proofs of other implications are almost equal to those in the proof of Theorem~\ref{PhiD}.
\par
Assume that $C_{\mathcal{ B}}(X)$ is countably dense in $\langle
\USCB{X},\tau_{\mathcal B}\rangle$, $\langle
\USCB{X},\tau_{\mathcal B}\rangle$ is separable and satisfies the
selection principle S${}_1(\tilde \Phi^{\uparrow}_{\mathcal
B},\Omega_{{\snula},{\mathcal B}})$. Thus, there exists a~
countable set $D=\{d_n:n\in\omega\}$ dense in $\langle
\USCB{X},\tau_{\mathcal B}\rangle$. Since $C_{\mathcal{ B}}(X)$ is
countably dense in $\langle \USCB{X},\tau_{\mathcal B}\rangle$,
for every $n\in\omega$ there exists a~countable set
$D_n=\{d_{n,m}:m\in\omega\}\subseteq C_{\mathcal{ B}}(X)$ such
that $d_n\in \overline{D_n}$ for each $n\in\omega$.
\par
Let $\{S_{n,m,l}:n,m,l\in\omega\}$ be a~sequence of subsets of
$\langle \USCB{X},\tau_{\mathcal B}\rangle$, each $S_{n,m,l}$
being in $\tilde\Phi^{\uparrow}_{\mathcal B}$. We can apply the
sequence selection principle
S${}_1(\tilde\Phi^{\uparrow}_{\mathcal
B},\Omega_{{\snula},{\mathcal B}})$ to the sequence
\[\seq{\{h-d_{n,m}:h\in S_{n,m,l}\}}{l}.\]
For every $l\in\omega$ we obtain $d_{n,m,l}\in S_{n,m,l}$ such that
\[\constant{0}\in\overline{\{d_{n,m,l}-d_{n,m}:l\in\omega\}}.\]
Then
\[d_{n,m}\in\overline{\{d_{n,m,l}:l\in\omega\}}.\]
Thus  $\{d_{n,m,l}:n,m,l\in\omega\}$ is the desired  countable dense set.
\qed
\par

\medskip

Since no infinite Hausdorff topological space has the covering
property S${}_1({\mathcal O},\Omega_{\mathcal B})$, we obtain
\begin{theorem}\label{noAD}
$\langle \USCB{X},\tau_{\mathcal B}\rangle$ does not have the
property {\rm S}${}_1({\mathcal P},{\mathcal D})$ for any
topological space $X$.
\end{theorem}

We generalize the main results of
\cite{BL1,BO,Osi,Osi1,Osi3,Osi4,Osi6}.

\begin{theorem}
Assume that $\Phi$ is one of the symbols $\Omega$ and $\Gamma$,
and $\Psi$~is one of the symbols ${\mathcal O}$, $\Omega$,
$\Gamma$. Then for any couple $\langle \Phi_{\mathcal
B},\Psi_{\mathcal B}\rangle$ different from
$\langle\Omega_{\mathcal B},{\mathcal O}\rangle$, a~topological
space $X$ is an~{\rm S}$_1(\Phi_{\mathcal B},\Psi_{\mathcal
B})$-space if and only if $\langle \USCB{X},\tau_{\mathcal
B}\rangle$ satisfies the selection principle
S${}_1(\Phi_{{\snula},{\mathcal B}},\Psi_{{\snula},{\mathcal
B}})$.
\par
If $X$ has the property ($\varepsilon_{\mathcal B}$) then the
equivalence holds true for the couple $\langle \Omega_{\mathcal
B},{\mathcal O}\rangle$ as well.
\end{theorem}

\medskip

\begin{corollary}\label{Cor6}
Assume that $\Phi$ is one of the symbols $\Omega$ and $\Gamma$,
and $\Psi$~is one of the symbols ${\mathcal O}$, $\Omega$,
$\Gamma$. Then for any couple $\langle \Phi_{\mathcal
B},\Psi_{\mathcal B}\rangle$ different from
$\langle\Omega_{\mathcal B},{\mathcal O}\rangle$, a~topological
space $X$ is an~{\rm S}$_1(\Phi_{\mathcal B},\Psi_{\mathcal
B})$-space if and only if for every $h\in\USCB{X}$ the family
$\langle \USCB{X},\tau_{\mathcal B}\rangle$ satisfies the
selection principle S${}_1(\Phi_{h,{\mathcal B}},\Psi_{h,{\mathcal
B}})$.
\end{corollary}

\section{Sequentially dense selectors of $\langle USC^*(X),\tau_{\mathcal{B}}\rangle$ and $\langle USC^*(X),\tau_{\mathcal{B}^s}\rangle$}
\par

\medskip
\begin{theorem}\label{PhiS}
Let $\Phi=\Omega$ or $\Phi=\Gamma$. Then for any couple
$\langle\Phi,\Psi\rangle$ different from $\langle \Omega_{\mathcal
B},{\mathcal O}\rangle$, the following are equivalent:
 \begin{enumerate}
\item[{\rm a)}] $\langle \USCB{X},\tau_{\mathcal B}\rangle$
satisfies the selection principle {\rm S}${}_1(\tilde
\Phi^{\uparrow}_{\mathcal B},{\mathcal S})$. \item[{\rm b)}]
$\langle \USCB{X},\tau_{\mathcal B}\rangle$ is sequentially
separable and the topological space~$X$ possesses the covering
property {\rm S}${}_1(\Phi_{\mathcal B},\Gamma_{\mathcal B})$.
 \item[{\rm c)}]  $\langle \USCB{X},\tau_{\mathcal B}\rangle$ is sequentially separable and satisfies the selection
principle {\rm S}${}_1(\Phi_{h,{\mathcal B}},\Gamma_{h,{\mathcal
B}})$ for every $h\in\USCB{X}$. \item[{\rm d)}] $\langle
\USCB{X},\tau_{\mathcal B}\rangle$ is sequentially separable and
satisfies the selection principles {\rm
S}${}_1(\Gamma_{{\snula},{\mathcal B}},\Gamma_{{\snula},{\mathcal
B}})$ and {\rm S}${}_1(\tilde \Phi^{\uparrow}_{\mathcal
B},\Gamma_{h,{\mathcal B}})$ for every $h\in\USCB{X}$.
\end{enumerate}
\end{theorem}
\pf
\newline
${\rm a)}\to {\rm b)}$. Let $\{{\mathcal
U}_n:n\in\omega\}\subseteq \Phi_{\mathcal B}$. We may assume that
${\mathcal U}_{n+1}$ is a~refinement of ${\mathcal U}_n$ for each
$n\in\omega$. If $\Phi=\Gamma$ we may also assume that  for every
$n\in\omega$, the cover ${\mathcal U}_n=\{U_i^n:i\in\omega\}$ is
a~countable family.
\par
We define the sets $S_n$ by \eqref{Sn}. By Lemma \ref{lem1},
$S_n\in\tilde\Phi^{\uparrow}_{\mathcal B}(\langle
\USCB{X},\tau_{\mathcal B}\rangle)$. We apply the selection
principle {\rm S}${}_1(\tilde\Phi^{\uparrow}_{\mathcal
B},{\mathcal S})$ and for every $n$ we obtain
a~function~$f_{U_n,h_n}\in S_n$ such that
$\{f_{U_n,h_n}:n\in\omega\}$ is sequentially dense in $\langle
\USCB{X},\tau_{\mathcal B}\rangle$. For every $n$ we shall find
a~set $V_n\in{\mathcal U}_n$ such that $\{V_n:n\in\omega\}$ is a
${\mathcal B}$-$\gamma$-cover.
\par
Evidently there exists an~increasing sequence $\seq{n_k}{k}$ such
that $f_{U_{n_k},h_{n_k}}\to \nula$. We set $V_{n_k}=U_{n_k}$. If
$n_k<n<n_{k+1}$, then by~(\ref{PhiPhih}) we can find a~set
$V_n\in{\mathcal U}_n$ such that $U_{n_{k+1}}\subseteq V_n$.
\par
Let $B\in {\mathcal B}$. Since $W=\{g\in\USCB{X}: \forall x\in B$
$|g(x)|<1 \}$ is a~neighborhood of $\nula$, there exists an~$k_0$
such that $f_{U_{n_k},h_{n_k}}\in W$ for each $k\geq k_0$. If
$g\in W$ then $g(x)<1$ for each $x\in B$. Thus for $k\geq k_0$ we
have $f_{U_{n_k},h_{n_k}}(x)<1$ for each $x\in B$. Therefore
$B\subseteq U_{n_k}=V_{n_k}$. By the choose of $V_n$ for
$n\notin\{n_k:k\in\omega\}$ we obtain that $B\subseteq V_n$ for
each $n\geq n_{k_0}$.
\par
The implication ${\rm b)}\to{\rm c)}$ follows by Corollary
\ref{Cor6}.
\par
The implication ${\rm c)}\to{\rm d)}$ is obvious by
\eqref{PhiPhih}.
\par
We prove the implication ${\rm d)}\to{\rm a)}$.
\par
Assume that there exists a~countable set
$\{d_n:n\in\omega\}\subseteq \USCB{X}$ sequentially dense in
$\langle \USCB{X},\tau_{\mathcal B}\rangle$ and $\langle
\USCB{X},\tau_{\mathcal B}\rangle$ satisfies the selection
principles S${}_1(\Gamma_{{\snula},{\mathcal
B}},\Gamma_{{\snula},{\mathcal B}})$ and  {\rm S}${}_1(\tilde
\Phi^{\uparrow}_{\mathcal B},\Gamma_{h,{\mathcal B}})$ for each
$h\in\USCB{X}$.
\par
Let $(S_{n,m}:n,m\in\omega)$ be a~sequence of subsets of $\langle
\USCB{X},\tau_{\mathcal B}\rangle$ all being in
$\tilde\Phi^{\uparrow}_{\mathcal B}$. For every $n$ we apply the
selection principle {\rm S}${}_1(\tilde \Phi^{\uparrow}_{\mathcal
B},\Gamma_{{d_n},{\mathcal B}})$ to the sequence
$(S_{n,m}:m\in\omega)$. Then for every $m\in\omega$ we obtain
$d_{n,m}\in S_{n,m}$, $d_{n,m}\geq d_n$, $d_{n,m}-d_n\in\USC{X}$
and such that $d_{n,m}-d_n\to \nula$ ($m\rightarrow \infty$).
\par
We show that the set  $\{d_{n,m}:n,m\in\omega\}$ is the desired
sequentially dense set.
\par
Indeed, if $h\in\USCB{X}$ then there exists an~increasing sequence
$\seq{n_k}{k}$ such that $d_{n_k}\to h$. Since for every $k$ we
have
\[d_{n_k,m}-d_{n_k}\to \nula,\]
by {\rm S}${}_1(\Gamma_{{\snula},{\mathcal
B}},\Gamma_{{\snula},{\mathcal B}})$ there exists a~sequence
$\seq{m_k}{k}$ such that
\[d_{n_k,m_k} -d_{n_k}\to \nula.\]
Thus \[d_{n_k,m_k}=(d_{n_k,m_k}- d_{n_k})+d_{n_k}\] converges to
$h$ ($k\rightarrow \infty$). \qed


\par

\medskip
Since no infinite Hausdorff topological space satisfies
S${}_1({\mathcal O},\Gamma_{\mathcal{B}})$ as above one can easily
prove
\begin{theorem}\label{noAS}
$\langle \USCB{X},\tau_{\mathcal{B}}\rangle$ does not have the
property {\rm S}${}_1({\mathcal P},{\mathcal S})$ for any
(infinite Hausdorff) topological space $X$.
\end{theorem}

\section{Fr\'echet-Urysohn property
in $\langle {\rm C}(X),\tau_{{\mathcal B}}\rangle$}


In \cite{BL21} and \cite{DCD}, remarkable studies of bornological
covering properties in $\langle {\rm C}(X),\tau_{{\mathcal
B}}\rangle$ were made. In this section we make some generalization
for the Fr\'echet-Urysohn property in $\langle {\rm
C}(X),\tau_{{\mathcal B}}\rangle$.

\par
If $\langle X,\upsilon\rangle$ is a~topological space, ${\mathcal
B}$ is a~bornology on $X$, we introduce the notion of
a~functionally separated-${\mathcal B}$-$\varphi$-cover for
$\varphi=\omega,\gamma$. A~cover~${\mathcal U}$ is a~{\emm
functionally separated-${\mathcal B}$-$\omega$-cover}, shortly
{\emm ${\mathcal B}^f$-$\omega$-cover},  if
\begin{eqnarray}\label{fs1}
{}&\hskip-3cm(\forall B\in{\mathcal B})(\exists U\in{\mathcal U})(\exists \delta_2>\delta_1>0)(\exists f\in\C{X})\\
{}&\hskip3cm((\forall  x\in B)\,f(x)<\delta_1\land (\forall x\in
X\setminus U)\,f(x)>\delta_2).\nonumber
\end{eqnarray}
One can easily see that the condition \eqref{fs1} is equivalent to
the condition
\begin{equation}\label{fs2}
(\forall B\in{\mathcal B})(\exists U\in{\mathcal U})(\exists
f\in\C{X})\,((\forall  x\in B)\,f(x)=0\land (\forall x\in
X\setminus U)\,f(x)=1).
\end{equation}
A~cover~${\mathcal U}$ is a~{\emm functionally
separated-${\mathcal B}$-$\gamma$-cover}, shortly {\emm ${\mathcal
B}^f$-$\gamma$-cover}, if ${\mathcal U}$ is infinite and for every
$B\in{\mathcal B}$ the set
\begin{equation}\label{fs3}
\{U\in{\mathcal U}:\lnot(\exists f\in\C{X})\,((\forall  x\in
B)\,f(x)=0\land (\forall x\in X\setminus U)\,f(x)=1)\}
\end{equation}
is finite. We denote by $\Omega_{{\mathcal B}^f}(X)$ and
$\Gamma_{{\mathcal B}^f}(X)$ the family of all open ${\mathcal
B}^f$-${\omega}$-covers and open ${\mathcal
B}^f$-${\gamma}$-covers of $X$, respectively. Then we have
\[\Gamma_{{\mathcal B}^f}(X)\subseteq\Omega_{{\mathcal B}^f}(X)\subseteq {\mathcal O}(X).\]
For a~Tychonoff space $X$ we obtain
\[\Gamma_{{\rm Fin}^f}(X)=\Gamma(X),\ \ \ \Omega_{{\rm Fin}^f}(X)=\Omega(X).\]
If $X$ is a~metric space then one can easily see that
\[\textnormal{every\ ${\mathcal B^s}$-$\varphi$-cover\ is\ a\ ${\mathcal B^f}$-$\varphi$-cover}.\]
For $X=\RR\setminus\QQ$ the opposite implication is false.
 \par
Similarly Theorem 12 in \cite{BL21} we obtain
\begin{lemma}\label{gn2}
Let ${\mathcal B}$ be a~bornology on a~topological space $\langle
X,\upsilon\rangle$ with a~closed base. The following are
equivalent:
\begin{enumerate}
\item[{\rm a)}] $X$ is a~${\rm S}_1(\Omega_{{\mathcal
B}^f},\Gamma_{{\mathcal B}^f})$ space. \item[{\rm b)}] Every open
${\mathcal B}^f$-$\omega$-cover has a countable ${\mathcal
B}^f$-$\gamma$-subcover.

If $X$ is a~metric space then similar equivalence holds true for
${\mathcal B}^s$-covers.

\item[{\rm e)}] every open ${\mathcal B}^f$-$\omega$-cover has a
countable ${\mathcal B}^f$-$\gamma$-subcover. \item[{\rm f)}] $X$
is a~${\rm S}_1(\Omega_{{\mathcal B}^f},\Gamma_{{\mathcal B}^f})$
space.
\end{enumerate}

\end{lemma}

Similarly Corollary 14 in \cite{BL21} we can prove
a~generalization of the Gerlits -- Nagy Theorem for the
topological space~$\langle {\rm C}(X),\tau_{{\mathcal B}}\rangle$.
\begin{theorem}\label{b2}
Let ${\mathcal B}$ be a~bornology with a~closed base on the
Tychonoff topological space $\langle X,\upsilon\rangle$. Then the
following are equivalent:
\begin{enumerate}
\item[{\rm a)}] $\langle {\rm C}(X),\tau_{{\mathcal B}}\rangle$ is
Fr\'echet-Urysohn. \item[{\rm b)}] every open ${\mathcal
B}^f$-$\omega$-cover has a countable ${\mathcal
B}^f$-$\gamma$-subcover. \item[{\rm c)}] $X$ is a~${\rm
S}_1(\Omega_{{\mathcal B}^f},\Gamma_{{\mathcal B}^f})$ space.
\item[{\rm d)}]  $\langle {\rm C}(X),\tau_{{\mathcal B}}\rangle$
is strictly Fr\'echet-Urysohn.
\end{enumerate}
\end{theorem}


\section*{Acknowledgments}
Professor Lev Bukovsk\'y was a wonderful man and a great
scientist. It was a great honor for me to work together. I hope
that his scientific ideas will be continued in further research by
mathematicians working in this field.

\end{document}